\documentclass[12pt]{amsart}

\usepackage[left=15mm, right=15mm, top=15mm, bottom=15mm]{geometry}

\usepackage{float}
\usepackage[font=small,labelsep=period, figurename={Рис.}]{caption}

\usepackage[english]{babel}
\usepackage[T1]{fontenc}
\usepackage[utf8]{inputenc}

\theoremstyle{plain}
    \newtheorem{theorem}{Theorem}[section]
    \newtheorem{lemma}[theorem]{Lemma}
    \newtheorem{proposition}[theorem]{Proposition}

\theoremstyle{definition}

    \newtheorem{example}[theorem]{Example}

\usepackage[
    draft = false,
    unicode = true,
    colorlinks = true,
    allcolors = blue,
    hyperfootnotes = true,
    citecolor = red
]{hyperref}
\usepackage{amsfonts,amsmath,amssymb,amscd,amsthm,url}
\usepackage[inline]{enumitem}
\usepackage{graphicx,epsf,afterpage, wrapfig}
\usepackage{soul}
\usepackage{multicol}
\usepackage{array}
\usepackage{braket}
\usepackage{epigraph}
\usepackage{rotating, floatflt}
\usepackage{xstring} 

\usepackage[
backend=biber,
style=numeric-comp,
sorting=none,
giveninits=true
]{biblatex}
\usepackage{csquotes}

\usepackage{xcolor}

\usepackage{ulem}

\usepackage{tikz}
\usetikzlibrary{positioning, calc, intersections, through, arrows, matrix, chains, math}
\usetikzlibrary{arrows.meta}
\usetikzlibrary{shadows}
\usetikzlibrary{graphs}
\usetikzlibrary{graphs.standard}
\usetikzlibrary{patterns}
\usetikzlibrary{decorations.pathreplacing,calligraphy,backgrounds}

\DeclareFieldFormat{usera}{\href{https://arxiv.org/abs/#1}{\texttt{arXiv:#1}}}

\renewbibmacro*{finentry}{\printfield{usera}\newunit\finentry}

\renewbibmacro{in:}{%
  \iffieldundef{usera}
    {\printtext{\bibstring{in}\intitlepunct}}
    {}%
}

\newcommand\abs[1]{\ensuremath{\left\lvert#1\right\rvert}}

\DeclareMathOperator{\cov}{cov}

\DeclareMathOperator{\add}{add}

\newcommand{\Acal}{\mathcal{A}}

\newcommand{\Fcal}{\mathcal{F}}

\newcommand{\Kcal}{\mathcal{K}}
\newcommand{\Lcal}{\mathcal{L}}
\newcommand{\Mcal}{\mathcal{M}}
\newcommand{\Ncal}{\mathcal{N}}

\newcommand{\defing}[1]{\textbf{\emph{\mathversion{bold}#1}}}

\newcommand{\R}{\ensuremath{\mathbb{R}}}

\newcommand{\N}{\ensuremath{\mathbb{N}}}

\renewcommand{\geq}{\geqslant}

\renewcommand{\leq}{\leqslant}

\newcounter{mcnt}

\newcounter{wordcnt}

\begin{document}

\title{Extension of kappa-additive measures}

\author{Sviatoslav V. Dzhenzher}

\begin{abstract}
    We establish an extension theorem for \(\kappa\)-additive probability measures, where \(\kappa\) is an arbitrary uncountable cardinal.
    Starting from a \(<\!\!\kappa\)-algebra, we define an outer measure using coverings of cardinality at most \(\kappa\) and characterise measurable sets by approximation in outer measure of symmetric differences. We prove that these sets form a \(\kappa\)-algebra carrying a complete \(\kappa\)-additive extension of the original measure. The extension is unique on the generated \(\kappa\)-algebra, and its completion is the minimal domain admitting a complete \(\kappa\)-additive extension. We also establish attainment of the covering infimum and show why the closure assumption on the initial algebra is essential. Examples include a strict extension of classical Lebesgue measure under the additional hypothesis \(\operatorname{add}\mathcal N=\aleph_1<\operatorname{cov}\mathcal N\).
\end{abstract}

\thanks{\hspace{-5mm}
S.\,V. Dzhenzher: sdjenjer@yandex.ru. orcid: 0009-0008-3513-4312
\\
Moscow Institute of Physics and Technology 141701, Institutskii lane, 9, Dolgoprudny, Russia}

\maketitle
\thispagestyle{empty}


\vspace{3mm}

\section{Introduction}

The theory of $\kappa$-additive measures is well studied \cite{Jech-p1-ch10, Fremlin-measure-vol5, Szpilrajn}.
Mostly, $\kappa$-additive measures defined on all subsets of a cardinal $\kappa$ have been studied, which leads to the area of measurable cardinals.
However, sometimes there arise measures that have no need to be defined on the entire powerset of unity.
For example, in \cite{Takesaki1979, Blecher-Weaver}, completely additive quantum states are described that induce completely additive measures on uncountable cardinal. In those particular cases, the induced measures are still defined on powersets, but these examples show that the theory of completely additive measures is interesting itself.

The classical extension theorem constructs a $\sigma$-additive measure on a generated \(\sigma\)-algebra from a $\sigma$-additive measure on an algebra of sets. A natural question is whether a similar construction can preserve additivity over uncountable families. Such forms of additivity are studied in set-theoretic measure theory \cite{Jech-p1-ch10, Fremlin-measure-vol5} and also arise in the theory of states on operator algebras. Complete additivity of a state on the projections of a von Neumann algebra characterises normality \cite{Takesaki1979}, whereas weaker forms of infinite additivity allow singular states and are related to measurable-cardinal properties \cite{Blecher-Weaver}. These connections motivate the study of extension procedures that preserve a prescribed degree of additivity.

For an uncountable cardinal \(\kappa\), the extension problem involves both the additivity of the measure and the closure properties of its domain. Throughout this paper, \(\kappa\)-additivity refers to disjoint families of cardinality at most \(\kappa\), with the additivity identity required whenever their union belongs to the domain of the measure. Under this convention, even complete additivity on an arbitrary algebra need not imply \(\kappa\)-subadditivity. Consequently, replacing countable coverings by coverings of cardinality at most \(\kappa\) may fail to preserve the original measure. An example on the first uncountable ordinal illustrates this obstruction.

We prove that closure of the initial algebra under unions of fewer than \(\kappa\) sets provides sufficient structure for an extension theorem. Given a \(\kappa\)-additive probability measure on such an algebra, we construct an outer measure from coverings of cardinality at most \(\kappa\). A set is declared measurable if it can be approximated by sets from the original algebra in the outer measure of their symmetric difference. The resulting class is a \(\kappa\)-algebra, and the restricted outer measure is a complete \(\kappa\)-additive extension of the original measure.
A further property of the construction is the attainment of the covering infimum for every subset of the underlying space. In the classical extension procedure starting from a \(\sigma\)-additive measure on an algebra, an optimal countable covering need not exist. Under our assumptions, however, every (even non-measurable) set admits a covering by at most \(\kappa\) members of the initial algebra whose total measure attains the outer measure of the set.

The examples demonstrate that the construction can enlarge a measure space beyond its ordinary completion. In particular, under some additional hypothesis, it yields a strict \(\aleph_1\)-additive extension of classical Lebesgue measure. The resulting domain still excludes ordinary Vitali sets.

The structure of the article is as follows.
Section~\ref{s:main} introduces the definitions and states the main results.
In particular, Theorem~\ref{t:extension} provides the result on the extension of a measure, analogous to that of Lebesgue;
Theorem~\ref{t:unique-caratheodory} and Theorem~\ref{t:unique-complete} provide the uniqueness of the obtained extensions.
After that, Sections~\ref{s:proof-ext}--\ref{s:proof-compl} prove the extension, uniqueness, and completion theorems.
Section~\ref{s:examples} presents examples showing the non-triviality of the obtained constructions.
Finally, Section~\ref{s:concl} discusses further directions.

\section{Background, main notions and results}\label{s:main}

For a set \(A\), by $2^A$ we denote its powerset, and by $\abs{A}$ we denote the cardinality of $A$.
We assume ZFC, so all cardinal numbers are comparable.
When we write \(A\subset B\), we mean \(A\subseteq B\).

Let \((\Omega,\Acal)\) be a \emph{measurable space}: that is, \(\Acal\subset 2^\Omega\) is an algebra of subsets of \(\Omega\in\Acal\).
Let \(m\colon\Acal\to[\,0,1\,]\) be a measure, which means a \emph{finitely additive} map; that is,
\[
    m(A\sqcup B) = m(A)+m(B)
    \quad\text{for any}\quad A,B\in\Acal.
\]
It is clear that the additivity extends to any finite disjoint unions.
We always presume that \(m(\Omega)=1\); we could consider any bounded measures \(\Acal\to[\,0,+\infty)\), but rescaling argument easily reduces the case to the probability measures \(\Acal\to[\,0,1\,]\).

Recall that for a family \(\{r_i\}_{i\in I}\) of non-negative reals,
\[
    \sum_{i\in I} r_i = \sup\left\{\sum_{i \in F} r_i\;:\; F\subset I,\,\abs{F} < \infty\right\}.
\]
The sum is finite only if the collection of strictly positive reals is at most countable.

Let \(\kappa\) be an uncountable cardinal.
Unlike that of \cite{Jech-p1-ch10}, but consistent with \cite{Blecher-Weaver}, we say that a measure \(m\colon\Acal\to[\,0,1\,]\) is said to be \defing{$\kappa$-additive} if
\[
    m\left(\bigsqcup_{i\in I}A_i\right) = \sum_{i\in I}m(A_i)
    \quad\text{for any}\quad \abs{I} \leq \kappa;
\]
in this definition, as in the definitions below, when we write \(m(A)\), we presume that this value is defined~--- in other words, that \(A\in\Acal\).
In \cite{Jech-p1-ch10}, $\kappa$-something meant \(\abs{I}<\kappa\); in these cases, following \cite{Blecher-Weaver}, we will write \emph{$<\!\!\kappa$-additivity}, and $<\!\!\kappa$-other things.
So, in our notation, we could write $\kappa$-additivity as $\leq\!\!\kappa$-additivity, but we do not want to.
Besides, $\kappa$-additivity is consistent with $\sigma$-additivity.

Analogously, $m$ is said to be \defing{$\kappa$-subadditive} if
\[
    m(A) \leq \sum_{i\in I}m(A_i)
    \quad\text{for any}\quad \abs{I} \leq \kappa,\quad
    A \subset \bigcup_{i\in I}A_i,
\]
and $m$ is said to be \defing{$\kappa$-superadditive} if
\[
    m(A) \geq \sum_{i\in I}m(A_i)
    \quad\text{for any}\quad \abs{I} \leq \kappa,\quad
    A \supset \bigsqcup_{i\in I}A_i.
\]
Finally, $m$ is said to be \defing{completely additive} if it is $\kappa$-additive for any cardinal $\kappa$.
It is clear that complete additivity is equivalent to $\abs{\Acal}$-additivity, or, generally, to $\abs{2^\Omega}$-additivity.
The \emph{complete subadditivity} and \emph{complete superadditivity} are defined analogously.

An algebra \(\Acal\) is called a \defing{$<\!\!\kappa$-algebra} if it is closed under unions of strictly fewer than $\kappa$ sets.
It is clear from De Morgan's laws that any $<\!\!\kappa$-algebra is closed under intersections of strictly fewer than $\kappa$-sets.
The notion of \defing{$\kappa$-algebra} is defined analogously.

The following result is well known for countable properties.

\begin{proposition}\label{p:add-sub-super}
    Let \(m\colon\Acal\to[\,0,1\,]\) be a (finitely additive) measure of an algebra~$\Acal$.
    Then
    \begin{enumerate}[label=(\roman*), ref=\roman*, nosep]
        \item\label{e:c-super} $m$ is completely superadditive,
        \item\label{e:k-sub} $\kappa$-subadditivity implies $\kappa$-additivity for any $\kappa$,
        \item\label{e:c-sub} complete subadditivity implies complete additivity.
        \item\label{e:k-equiv} If, additionally, \(\Acal\) is a $<\!\!\kappa$-algebra, then $\kappa$-additivity and $\kappa$-subadditivity are equivalent.
    \end{enumerate}
\end{proposition}

\begin{proof}
    Property~\eqref{e:c-super} follows from the definition of uncountable sums, supremum, and the fact that any measure is finitely superadditive.
    Implications
    \(\eqref{e:c-super} \Rightarrow\eqref{e:k-sub} \Rightarrow \eqref{e:c-sub}\) are clear.    
    It remains to prove~\eqref{e:k-equiv}.
    That is, if $\Acal$ is a $<\!\!\kappa$-algebra, and $m$ is $\kappa$-additive, then $m$ is $\kappa$-subadditive.
    So, let \(A\subset \bigcup_{i\in I}A_i\) for \(\abs{I}\leq \kappa\).
    By replacing \(A_i\) with \(A\cap A_i\), the sum \(\sum_{i\in I}m(A_i)\) can only decrease.
    So, it is sufficient to consider \(A = \bigcup_{i\in I}A_i\).
    Well order $I$ so that its proper initial segments have cardinality less than $\kappa$.
    Then, putting
    \[
        B_i := A_i \setminus \bigcup_{j < i} A_j \in \Acal,
    \]
    we obtain \(A = \bigsqcup_{i\in I}B_i\), and
    \[
        m(A) = \sum_{i \in I} m(B_i) \leq \sum_{i\in I} m(A_i). \qedhere
    \]
\end{proof}

\begin{example}[Complete additivity may not imply $\aleph_1$-subadditivity]\label{ex:not-sub}
    Let \(\Omega:=\omega_1\) be the first uncountable ordinal, which is also the well-ordered set of all at most countable ordinals.
    For \(\alpha<\omega_1\), denote by
    \[
        [\,0,\alpha\,] := \{\beta < \omega_1:\beta\leq\alpha\}
    \]
    the initial segment of $\omega_1$.
    Let $\Acal$ be the algebra generated by all \([\,0,\alpha\,]\).
    It is clear that \(\Acal = \Acal_{fin} \sqcup \Acal_{cofin}\), where the ring \(\Acal_{fin}\) consists of all subsets of the form
    \[
        \bigsqcup_{k=1}^n (\alpha_k, \beta_k\,]\quad\text{or}\quad
        [\,0,\beta_0\,] \sqcup \bigsqcup_{k=1}^n (\alpha_k, \beta_k\,],
        \qquad \alpha_k<\beta_k < \omega_1,\; n\in\N;
    \]
    and the family \(\Acal_{cofin}\) consists of all subsets of the form
    \[
        A \sqcup (\alpha,\omega_1),\qquad A\in\Acal_{fin}, \;\alpha<\omega_1.
    \]

    Let \(m\colon\Acal\to\{0,1\}\) be
    \[
        m(A) := \begin{cases}
            0, &A\in\Acal_{fin}, \\
            1, &A\in\Acal_{cofin}.
        \end{cases}
    \]
    Then $m$ is completely additive.
    Indeed, it is sufficient to prove that $m$ is $\aleph_1$-additive, but this follows from the finite additivity, since
    \[
        A = \bigsqcup_{i \in I} A_i
    \]
    for non-empty elements \(A_i\in\Acal\) implies that \(I\) is finite.

    But $m$ is not $\aleph_1$-subadditive, since
    \[
        \Omega = \omega_1 = \bigcup_{\alpha < \omega_1} [\,0,\alpha\,]
    \]
    and
    \[
        m(\Omega) = 1 > 0 = \sum_{\alpha < \omega_1} m[\,0,\alpha\,].
    \]
\end{example}

The previous example shows that we need to control the degree of the additivity of a measure.
So, below we will speak of $\kappa$-additivities for particular $\kappa$, but not of general complete additivities.

So, let $\kappa$ be an uncountable cardinal, and \(m\colon\Acal\to[\,0,1\,]\) be a $\kappa$-additive measure on an algebra $\Acal$.
We define the \defing{outer measure} \(\mu^*\colon 2^\Omega\to[\,0,1\,]\) by
\[
    \mu^*_\kappa(A) = \mu^*(A) := \inf\left\{\sum_{i\in I}m(A_i) \ :\ A\subset \bigcup_{i\in I}A_i, \; A_i\in\Acal, \; \abs{I}\leq\kappa\right\}.
\]
The motivation for such a definition is that we want to extend $m$ with the preservation of $\kappa$-additivity.
It is clear that \(\mu^*\) is bounded by unity, since $\Omega$ covers any set. It is also clear that the outer measure is monotone:
\[
    \mu^*(A) \leq \mu^*(B)
    \quad\text{for any}\quad A\subset B \subset\Omega.
\]

\begin{example}[Outer measure mail fail to extend the original measure on an algebra]\label{ex:outer-spoils}
    Take the algebra $\Acal$ and the measure $m$ from Example~\ref{ex:not-sub}.
    Since
    \[
        \Omega = \omega_1 = \bigcup_{\alpha < \omega_1} [\,0,\alpha\,]
    \]
    is covered by $\aleph_1$ subsets of zero measure, we have
    \[
        \mu^*_{\aleph_1}(\Omega) \leq \sum_{\alpha<\omega_1} m[\,0,\alpha\,] = 0 \neq m(\Omega).
    \]
    Thus, \(\mu^*_{\aleph_1} = 0\) for all subsets.
\end{example}

Example~\ref{ex:outer-spoils} shows that the analogues of the Carath\'eodory's extension theorem may not work even for a completely additive measure on an algebra.
The following lemma relies very heavily on the fact that \(\Acal\) is a \(<\!\!\kappa\)-algebra.

\begin{lemma}\label{l:outer}
    Let $\kappa$ be an uncountable cardinal.
    Let \(m\colon\Acal\to[\,0,1\,]\) be a $\kappa$-additive measure on a $<\!\!\kappa$-algebra $\Acal$.
    Let \(\mu^*=\mu^*_\kappa\) be the outer measure constructed by~$m$.
    Then
    \begin{enumerate}[label=(\roman*), ref=\roman*, nosep]
        \item\label{e:out-restr} the restriction \(\mu^*|_\Acal\) of \(\mu^*\) to $\Acal$ coincides with $m$; that is, $\mu^*$ extends $m$;
        \item\label{e:out-attain} the infimum in the definition of \(\mu^*\) is attained on some cover, and is equal to $m(P)$ for some $P\in \Acal$;
        \item\label{e:out-k-sub} $\mu^*_\kappa$ is $\kappa$-subadditive.
    \end{enumerate}
\end{lemma}

\begin{proof}
    Let us prove~\eqref{e:out-restr}.
    Let \(A \in \Acal\).
    Then \(\mu^*(A) \leq m(A)\) is clear, since $A$ covers $A$.
    On the other hand, for any cover \(A \subset \bigcup_{i \in I}A_i\) with \(\abs{I} \leq \kappa\), by Proposition~\ref{p:add-sub-super}.\ref{e:k-equiv}, $m$ is $\kappa$-subadditive, which means \(m(A) \leq \sum_{i\in I}m(A_i)\). This proves \(m(A) \leq \mu^*(A)\).

    Let us prove~\eqref{e:out-attain}.
    Let \(A\subset \Omega\).
    By the definition of infimum, for any \(\varepsilon > 0\), there exists a cover
    \[
        A \subset \bigcup_{i \in I_\varepsilon}A_{i,\varepsilon}
    \]
    such that
    \[
        \sum_{i\in I_\varepsilon} m(A_{i,\varepsilon}) < \mu^*(A) + \varepsilon.
    \]
    In particular, this sum is finite.
    Hence, the set \(I'_\varepsilon\subset I_\varepsilon\) of indices of sets of positive measure is at most countable.
    Denote their union by
    \[
        P_\varepsilon := \bigcup_{i \in I'_\varepsilon} A_{i,\varepsilon}.
    \]
    Since \(\Acal\) is a $\sigma$-algebra, \(P_\varepsilon \in \Acal\).
    Denote
    \[
        P := \bigcap_{\varepsilon \to 0} P_\varepsilon,
    \]
    where the intersection goes over a countable set of $\varepsilon$, and thus \(P \in \Acal\).
    Denote
    \[
        N := \bigcup_{\substack{\varepsilon \to 0 \\ i \in I_\varepsilon \setminus I'_\varepsilon}} A_{i,\varepsilon},
    \]
    where the union goes over the same $\varepsilon$ as in $P_\varepsilon$.
    This is the cover of cardinality at most \(\aleph_0\cdot \kappa = \kappa\).
    It is clear that
    \[
        A \subset P \cup N,
    \]
    From the definition of $P_\varepsilon$ and by $\kappa$-subadditivity of $m$ (again Proposition~\ref{p:add-sub-super}.\ref{e:k-equiv}), we have
    \[
        m(P_\varepsilon) \leq \sum_{i\in I_\varepsilon} m(A_{i,\varepsilon}) < \mu^*(A) + \varepsilon,
    \]
    and hence
    \[
        m(P) \leq \mu^*(A).
    \]
    On the other hand, since the cover \(P\cup N\) has the cardinality at most $\kappa$, and $N$ consists of sets of zero measure, by the definition of $\mu^*$ we have
    \[
        \mu^*(A) \leq m(P).
    \]
    Hence \(\mu^*(A) = m(P)\).

    Finally, let us prove~\eqref{e:out-k-sub}.
    Take any \(A\subset \bigcup_{i\in I}A_i\) for \(A,A_i\subset\Omega\) and \(\abs{I}\leq\kappa\).
    By~\eqref{e:out-attain}, there exist covers
    \[
        A_i \subset \bigcup_{j \in J_i} A_{i,j},
        \quad A_{i,j}\in\Acal, \; \abs{J_i}\leq\kappa,
    \]
    such that \(\mu^*(A_i) = \sum_{j\in J_i} m(A_{i,j})\).
    Since \(\kappa\cdot\kappa = \kappa\), we have
    \[
        \abs{\bigcup_{i\in I} J_i} \leq \kappa,
    \]
    and so
    \[
        \mu^*(A) \leq \sum_{i\in I}\sum_{j\in J_i} m(A_{i,j}) = \sum_{i\in I} \mu^*(A_i). \qedhere
    \]
\end{proof}

Thus, the picture is consistent with the countable case.
There, we had a \(\sigma\)-additive measure on an algebra, which is a \(<\!\!\aleph_0\)-algebra in our notation.
In our case, we have a \(\kappa\)-additive measure $m$ on a \(<\!\!\kappa\)-algebra.

Now we are ready to define measurable sets.
We will follow the intuitive definition analogously to \cite{BogachevSmolyanov-real-func, Pinelis-meas, Janssen}.
So, a set \(E \subset \Omega\) is said to be \defing{measurable} if for any \(\varepsilon>0\) there exists \(E_\varepsilon\in\Acal\) such that
\[
    \mu^*(E\bigtriangleup E_\varepsilon) < \varepsilon.
\]
We denote the class of measurable sets by \(\Mcal\), and the restriction of $\mu^*$ to \(\Mcal\) by \(\mu\).
It is clear that
\[
    \Mcal = \bigl\{E\subset \Omega\;:\; \inf_{A\in\Acal} \mu^*(E\bigtriangleup A)=0\bigr\} \supset \Acal.
\]

The following three theorems are classical for \(\kappa=\aleph_0(=\sigma)\), so we formulate them only for uncountable cardinals.

\begin{theorem}[Extension]\label{t:extension}
    Let $\kappa$ be an uncountable cardinal.
    Let \(m\colon\Acal\to[\,0,1\,]\) be a $\kappa$-additive measure on a $<\!\!\kappa$-algebra $\Acal$.
    Then the class \(\Mcal\) of measurable sets is a $\kappa$-algebra over $\Acal$, and \(\mu\) is a complete $\kappa$-additive extension of~$m$ to~$\Mcal$.
\end{theorem}

For an algebra \(\Acal\), denote by \(\kappa(\Acal)\) the minimal (by inclusion) $\kappa$-algebra over $\Acal$.
The correctness of the definition is clear: such the minimal $\kappa$-algebra provided by the facts that $2^\Omega$ is an appropriate $\kappa$-algebra, and that the intersection of all $\kappa$-algebras over $\Acal$ is a $\kappa$-algebra over $\Acal$.
For example, if \(\Acal\) is finite, then \(\kappa(\Acal) = \Acal\).
Another example: if \(\Acal\) is the algebra of finite unions of (possibly infinite) segments of $\R$, and \(\kappa = \R\), then \(\kappa(\Acal) = 2^\R\).

\begin{theorem}[Uniqueness]\label{t:unique-caratheodory}
    Let $\kappa$ be an uncountable cardinal.
    Let \(m\colon\Acal\to[\,0,1\,]\) be a $\kappa$-additive measure on a $<\!\!\kappa$-algebra $\Acal$.
    Then the restriction \(\mu|_{\kappa(\Acal)}\) of $\mu$ to \(\kappa(\Acal)\) provides the unique $\kappa$-additive extension of~$m$ to~$\kappa(\Acal)$.
\end{theorem}

\begin{theorem}[Completion]\label{t:unique-complete}
    Let $\kappa$ be an uncountable cardinal.
    Let \(m\colon\Acal\to[\,0,1\,]\) be a $\kappa$-additive measure on a $<\!\!\kappa$-algebra $\Acal$.
    Then the class $\Mcal$ of measurable sets is the minimal (by inclusion) $\kappa$-algebra over $\Acal$, on which a complete $\kappa$-additive extension of $m$ exists.
    Moreover, \(\mu\) is the unique such extension.
\end{theorem}

\section{Proof of Theorem~\ref{t:extension} (Extension)}\label{s:proof-ext}

\begin{lemma}\label{l:m-alg}
    Under the assumptions of Theorem~\ref{t:extension}, $\Mcal$ is an algebra.
\end{lemma}

\begin{proof}
    The idea is classical, so we provide it for the sake of completeness.
    
    Since \(\Mcal\) contains $\Acal$, we have \(\Omega \in \Mcal\).
    Let us show that \(\Mcal\) is an algebra. It is clear that it is closed under the complement, since
    \[
        E \bigtriangleup E_\varepsilon = (\Omega\setminus E) \bigtriangleup (\Omega\setminus E_\varepsilon).
    \]
    It remains to show that $\Mcal$ is closed under finite intersections.
    So, let \(A,B\in \Mcal\). Then, for any \(\varepsilon>0\),
    \[
        (A\cap B) \bigtriangleup (A_\varepsilon \cap B_\varepsilon) \subset (A\bigtriangleup A_\varepsilon) \cup (B \bigtriangleup B_\varepsilon),
    \]
    and, by the subadditivity of $\mu^*$ (Lemma~\ref{l:outer}.\ref{e:out-k-sub}),
    \[
        \mu^*\bigl((A\cap B) \bigtriangleup (A_\varepsilon \cap B_\varepsilon)\bigr) < 2\varepsilon. \qedhere
    \]
\end{proof}

\begin{lemma}\label{l:mu-f-add}
    Under the assumptions of Theorem~\ref{t:extension}, $\mu\colon\Mcal\to[\,0,1\,]$ is a finitely additive measure.
\end{lemma}

\begin{proof}
    Though the idea is classical, we do not know a reference for a simple intuitive proof, so we provide it here.
    Let \(A = B \sqcup C\) for some \(A,B,C\in\Mcal\).
    Since $\mu^*$ is subadditive (Lemma~\ref{l:outer}.\ref{e:out-k-sub}), it is sufficient to prove that
    \[
        \mu^*(A) \geq \mu^*(B) + \mu^*(C).
    \]
    For any \(\varepsilon>0\), take \(B_\varepsilon,C_\varepsilon \in \Acal\) by the definition of measurable sets.
    Then we may put
    \[
        A_{2\varepsilon} := B_\varepsilon \cup C_\varepsilon \in \Mcal.
    \]
    Indeed, we have already proved that $\Mcal$ is an algebra, so \(A_{2\varepsilon}\in\Mcal\), and it remains to use
    \[
        A \bigtriangleup A_{2\varepsilon} \subset (B\bigtriangleup B_\varepsilon) \cup (C\bigtriangleup C_\varepsilon)
    \]
    with the subadditivity of $\mu^*$ (Lemma~\ref{l:outer}.\ref{e:out-k-sub}).
    So, with the error $4\varepsilon$, it is sufficient to prove that
    \[
        \mu^*(A_{2\varepsilon}) \geq \mu^*(B_\varepsilon) + \mu^*(C_\varepsilon) - 2\varepsilon.
    \]
    But \(\mu^*\) restricted to $\Acal$ is the measure $m$.
    Then, by the inclusion-exclusion principle,
    \[
        \mu^*(A_{2\varepsilon}) = \mu^*(B_\varepsilon) + \mu^*(C_\varepsilon) - \mu^*(B_\varepsilon\cap C_\varepsilon).
    \]
    Since
    \[
        B_\varepsilon\cap C_\varepsilon \subset (B\bigtriangleup B_\varepsilon) \cup (C\bigtriangleup C_\varepsilon),
    \]
    we have, by the subadditivity of $\mu^*$ (Lemma~\ref{l:outer}.\ref{e:out-k-sub}), that \(\mu^*(B_\varepsilon\cap C_\varepsilon)<2\varepsilon\).
    Finally,
    \[
        \mu^*(A) \geq \mu^*(B) + \mu^*(C) - 6\varepsilon. \qedhere
    \]
\end{proof}

\begin{lemma}\label{l:m-lk-alg}
    Under the assumptions of Theorem~\ref{t:extension}, $\Mcal$ is a $<\!\!\kappa$-algebra.
\end{lemma}

\begin{proof}
    Let us show by the transfinite induction of \(\gamma < \kappa\) that \(U_\gamma := \bigcup_{\alpha < \gamma} E_\alpha \in \Mcal\) for \(E_\alpha \in \Mcal\).
    
    For a finite $\gamma$ it holds since \(\Mcal\) is an algebra.

    By the induction hypothesis, for \(\alpha < \gamma\) the unions
    \[
        U_\alpha = \bigcup_{\beta < \alpha} E_\beta
    \]
    lie in \(\Mcal\).
    Hence,
    \[
        D_\alpha := E_\alpha \setminus U_\alpha \in \Mcal,
    \]
    and thus
    \[
        U_\gamma = \bigcup_{\alpha < \gamma} E_\alpha = \bigsqcup_{\alpha < \gamma} D_\alpha.
    \]

    Now it becomes crucial that \(m\) is defined on an algebra.
    For any finite \(F \subset \gamma\), by the finite additivity of $\mu$ (Lemma~\ref{l:mu-f-add}), we have
    \[
        \sum_{i\in F} \mu(D_i) = \mu \left(\bigsqcup_{i\in F} D_i\right) \leq \mu(\Omega) = 1.
    \]
    Hence,
    \[
        \sum_{\alpha < \gamma} \mu(D_\alpha) = \sup_{F\subset\gamma}\sum_{i\in F} \mu(D_i) \leq 1.
    \]
    Then, by the definition of supremum, for any $\varepsilon>0$, there exists a finite subset \(F\subset\gamma\) such that
    \[
        \sum_{\substack{\alpha < \gamma \\ \alpha \notin F}} \mu(D_\alpha) < \varepsilon.
    \]
    Denote
    \[
        V := \bigsqcup_{\alpha \in F} D_\alpha \in \Mcal,
    \]
    where the inclusion holds since $\Mcal$ is an algebra, and the union is finite.
    Hence
    \[
        U_\gamma \bigtriangleup V = \bigsqcup_{\substack{\alpha < \gamma \\ \alpha \notin F}} D_\alpha,
    \]
    and, by $\kappa$-subadditivity of $\mu^*$ (Lemma~\ref{l:outer}.\ref{e:out-k-sub}),
    \[
        \mu^*(U_\gamma \bigtriangleup V) \leq \sum_{\substack{\alpha < \gamma \\ \alpha \notin F}} \mu(D_\alpha) < \varepsilon.
    \]
    But $V$ is measurable, hence, there exists \(V_\varepsilon\in\Acal\) such that
    \[
        \mu^*(V \bigtriangleup V_\varepsilon) < \varepsilon.
    \]
    Again by subadditivity of $\mu^*$ and from
    \[
        U_\gamma \bigtriangleup V_\varepsilon \subset
        (U_\gamma \bigtriangleup V) \cup (V \bigtriangleup V_\varepsilon)
    \]
    we obtain
    \[
        \mu^*(U_\gamma \bigtriangleup V_\varepsilon) < 2\varepsilon.
    \]
    The arbitrary choice of $\varepsilon$ proves the measurability of \(U_\gamma\).
\end{proof}

\begin{proof}[Proof of Theorem~\ref{t:extension}]
    The completeness of $\mu$ follows from the monotonicity of $\mu^*$ and from the fact that if \(\mu^*(E)=0\) then we always can take \(E_\varepsilon := \varnothing\).

    The $\kappa$-additivity of $\mu$ follows since
    \begin{itemize}[nosep]
        \item on a \(<\!\!\kappa\)-algebra, $\kappa$-additivity and $\kappa$-subadditivity of a finitely additive measure are equivalent (Proposition~\ref{p:add-sub-super}.\ref{e:k-equiv}),
        \item \(\mu^*\) is $\kappa$-subadditive (Lemma~\ref{l:outer}.\ref{e:out-k-sub}), and
        \item we have already proved that \(\Mcal\) is a $<\!\!\kappa$-algebra.
    \end{itemize}

    It remains to show that \(\Mcal\) is really a $\kappa$-algebra. The proof is just analogous to that one that $\Mcal$ is a $<\!\!\kappa$-algebra, with the only difference that \(\gamma=\kappa\), and that \(U_\alpha\) lie in $\Mcal$ not by the induction hypothesis, but since $\Mcal$ is a $<\!\!\kappa$-algebra.
\end{proof}

\section{Proof of Theorem~\ref{t:unique-caratheodory} (Uniqueness)}\label{s:proof-unique-k-alg}

\begin{proposition}\label{p:subadd}
    Let \(m\colon\Acal\to [\,0,1\,]\) be a subadditive map on an algebra~$\Acal$.
    Then for any \(A,B\in\Acal\),
    \[
        \abs{m(A)-m(B)} \leq m(A\bigtriangleup B).
    \]
\end{proposition}

\begin{proof}
    Holds since \(A \subset B \cup (A \bigtriangleup B)\).
\end{proof}

The author thanks Ivan Erlikh for suggesting (in the countable setting) the approximation argument used in the proof of Theorem~\ref{t:unique-caratheodory}.
See also \cite{Weizsacker-bmt} for this idea in case of $\sigma$-additive measures.

\begin{proof}[Proof of Theorem~\ref{t:unique-caratheodory}]
    Suppose that there is some other $\kappa$-additive extension \(\eta\) of \(m\) to \(\kappa(\Acal)\).
    Take any \(E \in \kappa(\Acal)\).
    Since \(E\) is measurable, for any \(\varepsilon>0\) there exists \(E_\varepsilon \in \Acal\) such that
    \[
        \mu(E\bigtriangleup E_\varepsilon) < \varepsilon.
    \]
    By Lemma~\ref{l:outer}.\ref{e:out-attain}, there exists some cover \(\bigcup_{i \in I} A_i \supset E\bigtriangleup E_\varepsilon\) with \(\abs{I}\leq\kappa\) and \(A_i\in\Acal\), such that
    \[
        \sum_{i \in I}m(A_i) = \mu(E\bigtriangleup E_\varepsilon).
    \]
    By $\kappa$-subadditivity of $\eta$ (Proposition~\ref{p:add-sub-super}.\ref{e:k-equiv}), we have
    \[
        \eta(E\bigtriangleup E_\varepsilon) \leq \sum_{i \in I}m(A_i) < \varepsilon.
    \]
    Now it remains to use the fact that \(\mu(E_\varepsilon) = m(E_\varepsilon) = \eta(E_\varepsilon)\) and apply Proposition~\ref{p:subadd} both to \(\mu\) and \(\eta\), in order to obtain
    \[
        \abs{\mu(E) - \eta(E)} \leq \abs{\mu(E) - \mu(E_\varepsilon)} + \abs{\eta(E_\varepsilon) - \eta(E)} \leq
        \mu(E\bigtriangleup E_\varepsilon) + \eta(E\bigtriangleup E_\varepsilon) < 2\varepsilon. \qedhere
    \]
\end{proof}

\section{Proof of Theorem~\ref{t:unique-complete} (Completion)}\label{s:proof-compl}

The proof is based on the following generalisation of the well known result describing the structure of measurable sets.

\begin{theorem}[Structure of measurable sets]\label{t:struct}
    Let $\kappa$ be an uncountable cardinal.
    Let \(m\colon\Acal\to[\,0,1\,]\) be a $\kappa$-additive measure on a $<\!\!\kappa$-algebra $\Acal$.
	Let \(\mu\) be the restriction of the outer measure on the $\kappa$-algebra $\Mcal$ of measurable sets.
    Then a set \(E \subset \Omega\) is measurable if and only if there exist some sets \(E' \subset E \subset E''\) from \(\kappa(\Acal)\) such that \(\mu(E') = \mu(E'')\).
\end{theorem}

\begin{proof}
	The ``if'' part follows by the completeness of $\mu$, since $E$ is the union of $E'$ and the subset \(E\setminus E'\) of the set \(E''\setminus E'\) of zero measure.
	So, we need to prove the ``only if'' part.
	Let $E \subset \Omega$ be a set (not necessarily measurable).
	By Lemma~\ref{l:outer}.\ref{e:out-attain}, there exists some cover \(E'':=\bigcup_{i\in I}A_i \supset E\) such that \(\abs{I} \leq \kappa\), \(A_i\in\Acal\), and
	\[
		\mu^*(E) = \sum_{i \in I} m(A_i) = \mu(E'').
	\]
	Analogously, we may find the cover \(\Omega \setminus E' \supset \Omega \setminus E\) such that
	\[
		\mu^*(\Omega \setminus E) = \mu(\Omega \setminus E').
	\]
	Hence, if $E$ is measurable, we obtain \(1 = \mu(E'') + 1 - \mu(E')\), which yields \(\mu(E') = \mu(E'')\).
\end{proof}

So, Theorem~\ref{t:struct} means that if we obtained a $\kappa$-additive extension \(\nu\) of $m$ to \(\kappa(\Acal)\), then \(\Mcal\) should be constructed by adding to $\kappa(\Acal)$ all subsets of sets of zero measure.
Formally, Theorem~\ref{t:unique-complete} follows from Theorem~\ref{t:struct} and Theorem~\ref{t:completion} below, applied to \(\nu = \mu|_{\kappa(\Acal)}\).

\begin{theorem}[Constructive completion]\label{t:completion}
	Let $\kappa$ be an uncountable cardinal.
    Let \(\nu\colon\Kcal\to[\,0,1\,]\) be a $\kappa$-additive measure on a $\kappa$-algebra $\Kcal$.
	Then the class \(\{E\bigtriangleup N\}\) of sets of the form \(E \bigtriangleup N\), where \(E \in \Kcal\) and $N$ is the subset of a set of zero measure, forms a $\kappa$-algebra,
	and the measure \(\overline{\nu}\) defined on this class by
	\[
		\overline\nu(E \bigtriangleup N) := \nu(E)
	\]
	is the $\kappa$-additive completion of $\nu$.
\end{theorem}

\begin{proof}
	First, let us provide the classical argument that three classes of sets of the form
	\[
		E \bigtriangleup N, \qquad E \cup N, \qquad E \sqcup N
	\]
	coincide.
	Both inclusions \(\{E \bigtriangleup N\} \supset \{E \sqcup N\}\) and \(\{E \cup N\} \supset \{E \sqcup N\}\) are clear.
	The inclusion \(\{E \cup N\} \subset \{E \sqcup N\}\) holds since \(E \cup N = E \sqcup (N \setminus E)\) and since \(N \setminus E \subset N\) is the subset of the same set of zero measure as $N$.
	The inclusion \(\{E \bigtriangleup N\} \subset \{E \sqcup N\}\) is proved as follows.
	Let \(N'' \supset N\) be a set of zero measure. Then
	\[
		E \bigtriangleup N = (E \setminus N'') \sqcup \bigl((E \cap (N'' \setminus N)) \sqcup N \setminus E\bigr),
	\]
	where \(E \setminus N'' \in \Kcal\) and \((E \cap (N'' \setminus N)) \sqcup N \setminus E \subset N''\).
	
	Now, let us show that \(\{E \bigtriangleup N\}\) is a $\kappa$-algebra.
	It is clear that it contains \(\Omega\).
	Below, by \(N''\in \Kcal\) we denote the set \(N''\supset N\) of zero measure.
	The closedness under symmetric difference is clear from definition.
	The closedness under finite intersections follows from the representation
	\[
		(E_1 \sqcup N_1) \cap (E_2 \sqcup N_2) = (E_1 \cap E_2) \sqcup (E_1 \cap N_2) \sqcup (N_1 \cap E_2) \sqcup (N_1 \cap N_2),
	\]
	where \(E_1 \cap E_2 \in \Kcal\), and the rest is the subset of the set \(N_1'' \cup N_2''\) of zero measure.
	Finally, the closedness under $\kappa$ unions follows from the representation
	\[
		\bigcup_{i\in I}(E_i \cup N_i) = \left(\bigcup_{i\in I}E_i\right) \cup \left(\bigcup_{i\in I}N_i\right),
	\]
	where the first union is from \(\Kcal\), and the second union is the subset of \(\bigcup_{i\in I}N_i''\).
	
	It remain to prove that \(\overline\nu\) is the $\kappa$-additive completion of $\nu$.
	The correctness of the definition follows since for \(E_1 \bigtriangleup N_1 = E_2 \bigtriangleup N_2\) we have \(E_1 \bigtriangleup E_2 = N_1 \bigtriangleup N_2\), and thus \(\nu(E_1) = \nu(E_2)\).
	Hence, $\overline\nu$ extends $\nu$.
	The completeness of \(\overline\nu\) holds since the class \(\{E\bigtriangleup N\}\) is constructed in this way.
	Finally, $\kappa$-additivity holds since for the family of $\kappa$ disjoint sets \(E_i \sqcup N_i\) we have
	\[
		\bigsqcup_{i\in I}(E_i \sqcup N_i) = \left(\bigsqcup_{i\in I}E_i\right) \sqcup \left(\bigsqcup_{i\in I}N_i\right),
	\]
	and hence
	\[
		\overline\nu\left(\bigsqcup_{i\in I}(E_i \sqcup N_i)\right) =
		\nu\left(\bigsqcup_{i\in I}E_i\right) = \sum_{i\in I} \nu(E_i) = 
		\sum_{i\in I} \overline\nu(E_i \sqcup N_i). \qedhere
	\]
\end{proof}

\section{Examples}\label{s:examples}

In this section, we provide three examples showing that Theorem~\ref{t:extension} gives non-trivial extension of the original measure, that is wider then the Lebesgue extension, and yet small enough not to cause issues with measurable cardinals.

\begin{example}
	Let \(\Acal\) be the $\sigma$-algebra of countable subsets of $\R$ and of their complements.
	Then \(\Acal\) is a \(<\!\!\aleph_1\)-algebra.
	Define \(m\colon\Acal\to[\,0,1\,]\) by
	\[
		m(A) := \begin{cases}
			0, &\text{$A$ is at most countable}, \\
			1, &\text{$\R\setminus A$ is at most countable}.
		\end{cases}
	\]
    Let us show that \(m\) is \(\aleph_1\)-additive if and only if the Continuum Hypothesis is false.
    Indeed, if the Continuum Hypothesis is true, then
    \[
        m(\R)=1 \neq 0 = \sum_{x\in\R}m\{x\},
    \]
    where the right sum is over \(\abs{\R} = \aleph_1\).
	Suppose that the Continuum Hypothesis is false.
    If \(\bigsqcup_{i \in I}A_i\) is at most countable, then all sets \(A_i\) are at most countable, and $\aleph_1$-additivity holds.
    On the other hand, if all the sets \(A_i\) are at most countable and \(\abs{I}\leq \aleph_1\), then
    \[
        \abs{\bigcup_{i\in I} A_i} \leq \aleph_0\cdot \aleph_1 = \aleph_1 < \abs{\R},
    \]
    which means that \(\bigcup_{i\in I} A_i\) must be at most countable in order to lie in \(\Acal\).

    Hence, we showed that \(m\) is an $\aleph_1$-additive measure on a \(<\!\!\aleph_1\)-algebra.
    Hence, by Theorem~\ref{t:extension}, there exists the extension \(\mu\) of $m$ to the $\aleph_1$-algebra $\Mcal$ of measurable sets.
    Moreover, in this case we can describe it explicitly.
    It is clear that
    \[
        \mu^*_{\aleph_1}(A) = \begin{cases}
            0, &\abs{A} \leq \aleph_1, \\
            1, &\abs{A} > \aleph_1.
        \end{cases}
    \]
    Hence, a set $E$ is measurable if and only if either \(\abs{E}\leq \aleph_1\) or \(\abs{\R\setminus E} \leq \aleph_1\).
    In other words, \(\aleph_1(\Acal) = \Mcal\) in this particular case.
\end{example}

\begin{example}
    In order to left the Continuum Hypothesis in peace, one should take in the previous example the space \(\omega_2\) instead of $\R$, and repeat all the constructions.
    Again, the $\aleph_1$-additive complete extension of $m$ will be obtained.
\end{example}

\begin{example}
    Let \(\lambda\) be the classical Lebesgue measure on \([\,0,1\,]\).
    Denote by \(\Lcal\) the class of Lebesgue measurable sets.
    Denote by
    \[
        \Ncal := \{N \in \Lcal\;:\; \lambda(N)=0\}
    \]
    the $\sigma$-ideal of sets of zero Lebesgue measure.
    Define the \emph{covering number}
    \[
        \cov \Ncal := \min\left\{\abs{\Fcal}\;:\; \Fcal\subset\Ncal, \; \bigcup \Fcal = [\,0,1\,]\right\},
    \]
    and define the \emph{additivity number}
    \[
        \add \Ncal := \min\left\{\abs{\Fcal}\;:\; \Fcal\subset\Ncal, \; \bigcup \Fcal \notin \Ncal\right\}.
    \]
    It is correct to write `minimum' here, since in ZFC cardinals are well ordered.
    Also, it is relatively consistent with ZFC that
    \[
        \add \Ncal = \aleph_1 < \cov\Ncal;
    \]
    see details in \cite{Fremlin-measure-vol5}, and also in \cite{Dzamonja-Plebanek}.
    
    In this paragraph we show that \(\aleph_1 < \cov\Ncal\) provides $\aleph_1$-additivity of $\lambda$.
    Let \(E_\alpha\) be Lebesgue measurable sets, and suppose that
    \[
        E := \bigsqcup_{\alpha<\omega_1} E_\alpha \subset [\,0,1\,]
    \]
    is also measurable.
    Then, by the complete superadditivity of $\lambda$ (Proposition~\ref{p:add-sub-super}.\ref{e:c-super}),
    \[
        \sum_{\alpha<\omega_1} \lambda(E_\alpha) \leq \lambda(E) \leq 1.
    \]
    Hence there exists at most countable family of $E_\alpha$ of positive measure.
    Denote their union by $P$.
    Then, due to the $\sigma$-additivity of $\lambda$,
    \[
        \lambda(P) = \sum_{\alpha<\omega_1} \lambda(E_\alpha).
    \]
    The set \(E\setminus P\) is covered by at most \(\aleph_1\) subsets of zero Lebesgue measure, hence
    \[
        \lambda(E\setminus P) = 0
    \]
    due to \(\aleph_1 < \cov\Ncal\).
    Finally,
    \[
        \lambda(E) = \lambda(P) + \lambda(E\setminus P) = \lambda(P) = \sum_{\alpha<\omega_1} \lambda(E_\alpha).
    \]

    In this paragraph we show that \(\add \Ncal = \aleph_1\) together with the $\aleph_1$-additivity of $\lambda$ provide that \(\Lcal\) is not an \(\aleph_1\)-algebra.
    Indeed, take \(N_\alpha \in \Ncal\) such that
    \[
        U := \bigcup_{\alpha<\omega_1} N_\alpha \notin\Ncal,
    \]
    which means \(\lambda^*(U)>0\). If \(\Lcal\) were an $\aleph_1$-algebra, then $U$ would be measurable, and $\lambda$ would be $\aleph_1$-subadditive (by Proposition~\ref{p:add-sub-super}.\ref{e:k-equiv}), and this would lead to \(\lambda(U) = 0\), which is the contradiction.

    Thus we showed that under the additional hypothesis
    \[
        \add \Ncal = \aleph_1 < \cov\Ncal,
    \]
    the classical Lebesgue measure is an $\aleph_1$-additive measure on a $<\!\!\aleph_1$-algebra.
    Hence we can \emph{extend the classical Lebesgue measure} to some complete $\aleph_1$-additive measure on an $\aleph_1$-algebra over $\Lcal$; cf.~\cite{Szpilrajn}.
    It is interesting that the classical arguments provide that the Vitali set will remain non-measurable even after the extension.
\end{example}

\section{Conclusion and further research}\label{s:concl}

We have established an extension procedure for \(\kappa\)-additive probability measures on \(<\!\!\kappa\)-algebras.
The construction uses coverings of cardinality at most \(\kappa\) and intuitive approximation by symmetric differences.
It produces a complete \(\kappa\)-additive measure, is unique on the generated \(\kappa\)-algebra, and yields the minimal domain admitting a complete extension of this kind. The examples show both the obstruction that occurs on arbitrary algebras and the possibility of obtaining extensions strictly larger than the ordinary completion. In particular, a suitable hypothesis on the Lebesgue null ideal gives a strict extension of classical Lebesgue measure without making every subset measurable.

A natural direction for further work is the study of integration with respect to these measures.
Blecher and Weaver \cite[Lemma 3.2]{Blecher-Weaver} use a transfinite version of the monotone convergence theorem to interchange integration and an uncountable sum of non-negative functions; see also \cite{Blecher-Labuschagne}.
Their argument concerns measures defined on entire power sets. The present construction provides measure spaces on proper \(\kappa\)-algebras in which the corresponding cardinal bound controls both measurable unions and unions of exceptional null sets. This suggests a systematic study of uncountable summation and convergence principles in this setting, together with their consequences for positive functionals. Further questions concern vector-valued and operator-valued analogues and their possible applications to states on operator algebras.

See some other applications of $\kappa$-additive measures in \cite{Dzhenzher26-qsing}, where quantum channels with singularising yet preserving $\sigma$-additivity properties are presented. The control of $\kappa$-additivity there would probably produce singularising quantum channels preserving $\kappa$-additivity of quantum states.
It is worth to mention that in \cite{Dzhenzher26-qsing}, measures were defined on entire powersets, in order to avoid the problems with measurability and integrability of functions. If we consider $\kappa$-additive measures on $\kappa$-algebras, then probably such problems of measurability can be avoided similarly to \cite{DzhenzherDzhenzherSakbaev26}.

\printbibliography

\end{document}